\documentclass{amsart}

\usepackage{amsmath,amssymb}
\usepackage{amsthm}
\usepackage[alphabetic]{amsrefs}
\usepackage{indentfirst}
\usepackage{mathrsfs}
\usepackage{graphicx}
\usepackage{hyperref}
\usepackage{xcolor}
\usepackage[margin=1.5in]{geometry}
\usepackage{graphicx,tikz}
\newtheorem{theorem}{Theorem}

\theoremstyle{definition}

\theoremstyle{remark}

\begin{document}

\title[]{A Remark on the Arcsine Distribution \\and the Hilbert transform}
\keywords{Hilbert transform, Arcsine distribution, Integral Operators.}
\subjclass[2010]{44A15, 60E99 (primary) and 33C45 (secondary).}

\author[]{Ronald R. Coifman}
\address{Department of Mathematics, Yale University, New Haven, CT 06511, USA}
\email{ronald.coifman@yale.edu, stefan.steinerberger@yale.edu}

\author[]{Stefan Steinerberger}

\thanks{C. R. is partially supported by grants from the NIH. S.S. is partially supported by the NSF (DMS-1763179) and the Alfred P. Sloan Foundation.}

\begin{abstract}
It is known that if $(p_n)_{n \in \mathbb{N}}$ is a sequence of orthogonal polynomials in $L^2([-1,1], w(x)dx)$, then the roots are distributed according to an arcsine distribution $\pi^{-1} (1-x^2)^{-1}dx$
for a wide variety of weights $w(x)$. We connect this to a result of the Hilbert transform due to Tricomi: if $f(x)(1-x^2)^{1/4} \in L^2(-1,1)$ and its Hilbert transform $Hf$ vanishes on $(-1,1)$, then the function $f$ is a multiple of the arcsine distribution
$$ f(x)  =  \frac{c}{\sqrt{1-x^2}}\chi_{(-1,1)} \qquad \mbox{where}~c~\in \mathbb{R}.$$
We also prove a localized Parseval-type identity that seems to be new: if $f(x)(1-x^2)^{1/4} \in L^2(-1,1)$ and $f(x) \sqrt{1-x^2}$ has mean value 0 on $(-1,1)$, then
$$  \int_{-1}^{1}{ (Hf)(x)^2 \sqrt{1-x^2} dx} =  \int_{-1}^{1}{ f(x)^2 \sqrt{1-x^2} dx}.$$
\end{abstract}

\maketitle

\section{Introduction }
\subsection{Introduction.} This short paper is concerned with properties of the Hilbert transform
$$ (Hf)(x) = \frac{1}{\pi}\int_{\mathbb{R}}{\frac{f(y)}{x-y} dy}$$
when interpreted as an operator acting on $L^2([-1,1])$ (this is sometimes called the 'finite' Hilbert transform $H_T$). One of the first people who
seems to have studied this problem is Tricomi \cite{tricomi} who was interested in solving the so-called airfoil equation $H_t f = g$. He observed
the following fundamental result.

\begin{theorem}[Tricomi \cite{tricomi}] Let $f(x)(1-x^2)^{1/4} \in L^2(-1,1)$. If the Hilbert transform $Hf$ vanishes identically on $(-1,1)$, then
$$ f(x)  =  \frac{c}{\sqrt{1-x^2}}\chi_{(-1,1)} \qquad \mbox{where}~c~\in \mathbb{R}.$$
\end{theorem}

We will give a particularly simple proof of Theorem 1 which has the advantage of also establishing a Plancherel-type theorem for $H_T$ that
seems to be new (see Theorem 2 in \S 1.3.). We also discuss a connection to orthogonal polynomials where Theorem 1 has a particularly nice interpretation and inverse problems for integral operators.

\subsection{Orthogonal polynomials.} It is well understood that, given a nonnegative weight on $(-1,1)$, the associated family
of orthogonal polynomials has roots whose distribution tend to the arcsine distribution (this dates back to Erd\H{o}s \& Turan \cite{erd}
in 1940, see also Erd\H{o}s \& Freud \cite{erd2}, Ullman \cite{ull} and Van Assche \cite{van}). Our main result provides a fairly natural
way to explain why that this is indeed the only natural smooth distribution that could have that property. If $p_n$ is a polynomial having
roots $\left\{x_1, \dots, x_n\right\} \subset (-1,1)$, then
$$ \frac{p_n'(x)}{p_n(x)} = \sum_{k=1}^{n}{\frac{1}{x-x_k}}$$
can be used to derive information about the distribution of the roots of $p_n'$. If the quantity was very big (in absolute values) in a region,
then it would repel any roots from locating there. This suggests that the only way to have the roots of $p_n'$ follow the same distribution
as the roots of $p_n$ is for the sum to somewhat cancel out; conversely, the sum is approximated by the Hilbert transform of the probability
measure of the distribution of roots which we would thus like to see vanish everywhere. It was this consideration, motivated by a recent dynamical interpretation
\cite{steini} of roots, that originally lead us to Theorem 1.\\

A specific recent application of Theorem 1 (partially inspired by this interpretation) in this context is as follows: suppose $p_n(x)$ is a polynomial of degree $n$ on the real line having only real roots.
Suppose furthermore that the roots are approximately distributed according to a smooth (possibly compactly supported) probability distribution $u_0(x)$. 
What can one say about the distribution of roots of $p_n'$?  Rolle's theorem implies that $p_n'$ has $n-1$ roots on the real line and that these roots interlace
with the roots of $p_n$. In particular, we can expect, as $n \rightarrow \infty$, that the roots of $p_n'$ are also distributed according to $u_0(x)$ and, more
generally, the same is true for the $k-$th derivative $p_n^{(k)}$. However, this is no longer true when $k$ is allowed to move with $n$: what can be said about
the distribution of roots of $p_n^{(0.1n)}$ or, more generally the distribution  $u(t,x)$ of $p_n^{(tn)}$ for $0 < t < 1$ depending on $u_0(x)$? The second author \cite{stein2} recently proposed a nonlinear, nonlocal transport
equation
$$ \frac{\partial u}{\partial t}  + \frac{1}{\pi} \frac{\partial}{\partial x}\left( \arctan{ \left( \frac{Hu}{ u}\right)} \right) = 0,$$
where $H$ is the Hilbert transform and the equation acts on $\mbox{supp}(u(t,x))$. The derivation is not rigorous but recovers the correct results for orthogonal polynomials
on $(-1,1)$, the family of Hermite polynomials (where the equation turns into a one-parameter family of shrinking semicircle distributions) and the family of Laguerre polynomials (where the equation turns into a one-parameter flow within the family of Marchenko-Pastur distribution).
Theorem 1 becomes significant in the first case since $u_t = 0$ locally requires $Hu = 0$ leading naturally to the arcsine distribution. A naturally variant of the equation on the torus has since
been studied by Granero-Belinchon \cite{granero} who established well-posed in certain spaces.

\subsection{Integral Operators.} There is a secondary motivation: while upper bounds on integral operators are well understood,
there is no such corresponding theory for lower bounds. One simple question one could ask is the following: let $f \in C^{\infty}_{c}(-1,1)$,
how big does the Hilbert transform have to be on, say, the interval $(2,3)$? A sharp result was given by Alaifari, Pierce and the
second author in \cite{al} (see also R\"uland \cite{angkana}) and reads
$$  \|H f\|_{L^2(2,3)} \geq c_1 \exp{\left(-c_2\frac{ \|f'\|_{L^2(-1,1)}}{\|f\|_{L^2(-1,1)}}\right)} \| f \|_{L^2(-1,1)},$$
The proof is far from stable (in the sense that it is not clear how to establish
arguments of this type for more general integral operators). 
The problem seems to be completely
open for general integral operators (we refer to \cite{roy} for sharp results for the Laplace transform
and the Fourier transform and to \cite{brett} for a dyadic model). 
When dealing with the Hilbert transform, the identity $\|Hf\|_{L^2(\mathbb{R})} = \|f\|_{L^2(\mathbb{R})}$ suggests a rephrasing of
the question: how does the Hilbert transform move the $L^2-$mass of a function around? Our main result in that direction reads as follows.
\begin{theorem} Let $f(x)(1-x^2)^{1/4} \in L^2(-1,1)$. If $f(x)\sqrt{1-x^2}$ has mean value 0, then
$$  \int_{-1}^{1}{ (Hf)(x)^2 \sqrt{1-x^2} dx} =  \int_{-1}^{1}{ f(x)^2 \sqrt{1-x^2} dx}.$$
\end{theorem}
This is a nice addition to the classical global $L^2-$isometry $\|Hf\|_{L^2(\mathbb{R})} = \|f\|_{L^2(\mathbb{R})}$ for compactly supported functions. 
We were surprised to not find this result in the literature, it seems so nice that one would expect that it should be known and not difficult to find.
Are there similar statements
(perhaps not identities but maybe inequalities) for more general singular integral operators of convolution type?
We refer to papers of Astala, P\"aiv\"arinta and Saksman \cite{finnish}, Bertola, Katsevich and Tovbis \cite{bert},  Katsevich \cite{kat}
Okada \& Elliott \cite{okada, okada1, okada2} and references therein for an overview of existing results.

\section{Proofs}

\subsection{A real-variable approach to Theorem 1 and Theorem 2.}
Before giving a rigorous proof in \S 2.2, we sketch an argument that is not entirely rigorous
(we do not specify the regularity of the function and freely interchange summation and integration).
However, the argument does motivate the rigorous proof in the sense that it can be interpreted
as the 'correct' argument (happening on the upper half disk) projected down onto the interval 
(this also explains the origin of several of the underlying identities).

\begin{proof}[Sketch of an argument]
 We introduce the function
$ f(x) \sqrt{1-x^2}$
which reduces the problem to showing that if the Hilbert transform of that function vanishes on $(-1,1)$, then that function is constant. We expand it into Chebychev polynomials
$$ f(x) \sqrt{1-x^2} = \sum_{k=0}^{\infty}{a_k T_k(x)}$$
and write
\begin{align*}
\frac{1}{\pi}\int_{-1}^{1}{\frac{f(y)}{x-y} dy} &= \frac{1}{\pi}\int_{-1}^{1}{\frac{f(y) \sqrt{1-y^2}}{(x-y) \sqrt{1-y^2}} dy} \\
&= \frac{1}{\pi}\int_{-1}^{1}{\sum_{k=0}^{\infty}{\frac{a_k T_k(y)}{(x-y) \sqrt{1-y^2}} dy}}.
\end{align*}
We want to show that this quantity vanishing implies that $f(x) \sqrt{1-x^2} = a_0$. 
The key ingredient is the identity
$$   \frac{1}{\pi}\int_{-1}^{1}{\frac{a_k T_k(y)}{(x-y) \sqrt{1-y^2}} dy} = a_k U_{k-1}(x)$$
where $U_{\ell}$ denotes the Chebychev polynomials of the second kind given by 
$$U_{0}(x) =1, U_{1}(x) = 2x \quad \mbox{and} \quad U_{\ell+1}(x) = 2x U_{\ell}(x) - U_{\ell-1}(x).$$
Writing
$$ f(x) \sqrt{1-x^2} = \sum_{k=0}^{\infty}{a_k T_k(x)},$$
we have that
\begin{align*}
\int_{-1}^{1}{ f(x)^2 \sqrt{1-x^2} dx} &= \int_{-1}^{1}{ \left(f(x) \sqrt{1-x^2}\right)^2 \frac{dx}{\sqrt{1-x^2}}} \\
&=\int_{-1}^{1}{ \sum_{k, \ell = 0}^{\infty}{ a_k T_k(x) a_{\ell} T_{\ell}(x)}   \frac{dx}{\sqrt{1-x^2}} } \\
&=  \sum_{k, \ell = 0}^{\infty}{  \int_{-1}^{1}{a_k T_k(x) a_{\ell} T_{\ell}(x)  }   \frac{dx}{\sqrt{1-x^2}}}\\
&= \sum_{k= 0}^{\infty}{  \int_{-1}^{1}{a_k^2 T_k(x)^2 }  \frac{dx}{\sqrt{1-x^2}}} \\
&= \frac{\pi}{2} \sum_{k=0}^{\infty}{ a_k^2}
\end{align*}
and the identity above shows that
\begin{align*}
 \int_{-1}^{1}{ (Hf)(x)^2 \sqrt{1-x^2} dx} &=  \int_{-1}^{1}{ \left( \sum_{k=1}^{\infty}{a_{k} U_{k-1}(x)} \right)^2  \sqrt{1-x^2} dx} \\
&=  \int_{-1}^{1}{  \sum_{k, \ell=1}^{\infty}{a_{k} U_{k-1}(x) a_{\ell} U_{\ell-1}(x) \sqrt{1-x^2} dx} }\\
&=  \sum_{k=1}^{\infty}{ \int_{-1}^{1}{  a_{k}^2 U_{k-1}(x)^2 \sqrt{1-x^2} dx} }\\
&= \frac{\pi}{2} \sum_{k=1}^{\infty}{ a_k^2}
\end{align*}
resulting in the desired identity
$$ \int_{-1}^{1}{ f(x)^2 \sqrt{1-x^2} dx} = \int_{-1}^{1}{ (Hf)(x)^2 \sqrt{1-x^2} dx} $$
whenever $f(x) \sqrt{1-x^2}$ has mean value 0. 
\end{proof}

The sketch of the argument has perhaps a somewhat uneasy feel to it: it is not a priori clear that integration and summation can
be exchanged like that and the algebraic identities seem to be coming out of nowhere: in truth, we are observing a fairly natural proof
in a different space projected down on the unit interval. The next section contains a more illuminating proof
(which, simultaneously, explains the origin of the identities used above).

\subsection{Proof of Theorem 1 and Theorem 2}
The main idea is
a substitution that can be found in a paper of Coifman \& Weiss \cite{coifman} relating the Hilbert transform on $(-1,1)$
to the conjugate function on the boundary of the disk (the trick is surely older than that); as such, the proof is algebraic in nature and unlikely to generalize to other integral
operators.
\begin{proof} Let $f(x)(1-x^2)^{1/4} \in L^2(-1,1)$ be given. We define a function $g:[-\pi, \pi]\rightarrow \mathbb{R}$ via
$$ g(\theta) = f(\cos{\theta}) \sqrt{1-(\cos{\theta})^2}.$$
We want to make sure that $g \in L^2(-\pi, \pi)$, a simple substitution $x = \cos{\psi}$ shows that
$$   \int_{-\pi}^{\pi}{ g(\theta)^2 d\theta} =   \int_{-\pi}^{\pi}{ f(\cos{\theta})^2 (1- (\cos{\theta})^2) d\theta} = \int_{-1}^{1}{ f(x)^2 \sqrt{1-x^2}dx} $$
which is finite by assumption. The, however, a simple substitution (carried out in \cite{coifman}) shows that with $x = \cos{\psi}$
$$ \mbox{p.v.}  \frac{1}{\pi}\int_{-1}^{1}{\frac{f(y)}{x-y} dy} = - \frac{1}{2 \sin{\psi}} \mbox{p.v.} \frac{1}{\pi} \int_{-\pi}^{\pi}{  g(\theta) \cot{ \left(\frac{\psi - \theta}{2} \right)} d\theta}.$$
However, that last principial value is merely the formula for the conjugate function of $g$. If $g$ has mean value 0 (this being
the orthogonality to the arcsine distribution),
then the conjugate function has the same $L^2-$norm as the function $g$. This implies the result.
\end{proof}

\textbf{Acknowledgment.} We are grateful to the referee for refering us to the work of Tricomi and several helpful suggestions.

\end{document}